\documentstyle[amscd,amssymb,verbatim,12pt,righttag]{amsart}
 \numberwithin{equation}{section}
 \renewcommand{\rm}{\normalshape}
 
\theoremstyle{plain}
\newtheorem{thm}{Theorem}[section]
\newtheorem{lemma}[thm]{Lemma}
\newtheorem{pro}[thm]{Proposition}
\newtheorem{cor}[thm]{Corollary}

\def\prf{{\bf Proof.\ }} 

\begin{document}
\title{ The Pressure Function for Products\\
 of Non-negative Matrices
}

\thanks{
\mbox{}\ \ \ The first author was partially supported by
a HK RGC grant in Hong Kong and the Special Funds for Major State
Basic Research Projects in China. The second author is supported by a HK RGC
grant.
}

\author{De-Jun FENG and Ka-Sing LAU}

\date {}


\maketitle
{\small {\bf Abstract.} Let $(\Sigma_A, \sigma)$ be a subshift of finite type and let $M(x)$ be a continuous function on $\Sigma_A$ taking values in the set of non-negative matrices. We extend the classical scalar pressure function to this new setting and prove the existence of the Gibbs measure and the differentiability of the pressure function. We are especially interested on the case where $M(x)$ takes finite values $M_1, \cdots, M_m$.  The pressure function reduces to $P(q):=\lim_{n\rightarrow \infty}\frac{1}{n} \log \sum_{J \in \sum_{A, n}} \|M_J\|^q$.  The expression is  important when we consider the multifractal formalism for certain iterated function systems with overlaps.}

\maketitle

\section{\bf Introduction}

Let $\sigma$ be the shift map on $\Sigma=\{1,2,\cdots, m\}^{\Bbb N}$, $m\geq 2$. As usual $\Sigma$ is endowed with the  metric 
$d(x,y)=m^{-n}$ where $x=(x_k),\ y=(y_k)$ and  $n$ is the smallest  of the $k$ such that $x_k \not = y_k$. Given an $m\times m$ matrix $A$ with entries $0$ or $1$, we consider the subshift of finite type $(\Sigma_A,\sigma)$ (see \cite{B}). We shall always assume that $A$ is primitive.

Suppose $M$ is a continuous function on $\Sigma_A$ taking values in the set of all non-negative $d\times d$ matrices. 
For $q\in {\Bbb R}$, we define the pressure function $P(q)$ of $M$ by 
\begin{equation}\label{n0}
P(q)=\lim_{n\rightarrow \infty} \frac{1}{n}\log 
\sum_{J\in \Sigma_{A,n}}\sup_{x\in [J]} \| M(x)M(\sigma x)
         \ldots M(\sigma^{n-1}x)\|^q,
\end{equation}
where $\Sigma_{A,n}$ denotes the set of all admissible indices of length $n$
 over $\{1,\ldots,m\}$; for 
 $ J =j_1 \cdots j_n \in\Sigma_{A,n}$, $[J]$ denotes the cylinder
 set $\{x=(x_i)\in \Sigma_A: \ x_i = j_i,\ 1\leq i\leq n\}$,
$\|\cdot\|$ denotes the matrix norm defined by $\|B\|:={\bf 1}^tB{\bf 1}$, ${\bf 1}^t=(1,1,\ldots,1)$. By using a subadditive argument, it is easy to show that for $q>0$, the limit in the above definition exists. With some additional conditions on the matrices (see Theorem 1.1), the limit exists for $q \in {\Bbb R}$.  

\bigskip

The pressure function of the scalar case (i.e., $M(x) = e^{\phi (x)}$ where $\phi(x)$ is a real valued function called the {\it potential} of the subshift) has been studied in great detail in statistical mechanics  and dynamical systems in conjunction with the Gibbs measure, the entropy and the variational principle (c.f., e.g., [B], [P], [R]); it has also been used to study the multifractal structure of the self-similar (or self-conformal) measures generated by iterated function systems (IFS) with no overlap (the {\it open set condition}) ([MU], [FL]). By identifying with the symbolic space, such  self-similar measure $\mu$ is actually a Gibbs measure and the pressure function is directly related to the scaling spectrum of $\mu$ [FL, Theorem 3.3].  In all the above cases, the pressure functions under consideration are differentiable (actually real analytic). This property is essential to investigate the phase transition in thermodynamics and for the  multifractal formalism in the dimension theory of fractals. 

\medskip

In the recent investigation of the self-similar measures generated by iterated function systems with {\it overlaps}, it is seen that in many interesting cases, such measure $\mu$ can be put into a vector form with a new non-overlapping IFS and with matrix weights ([LN1,2], [LNR], [Fe], [FeO]).  In this way the validity of the multifractal formalism depends on the differentiability of the pressure function $P(q)$ in (1.1) (more precisely (1.4) in the following) [LN2]. In another direction, the expression of the matrix product in (1.1) also appears in the study of the scaling functions in wavelet theory (the matrices are allowed to have negative entries) in the form of  $L^q$-joint spectral radius and the $L^q$-Lipschitz exponent ([DL1, 2], [LM]);  the problem of differentiability of the $P(q)$ also appears there.  So far there is no general theorem to guarantee this fact other than some special cases (e.g.,  [LN1], [FLN], [Fe], [FeO], [DL2]). 

\medskip

The main purpose of this paper is to consider the pressure functions and the Gibbs measures for the products of matrices. We first study the case that the matrices $M(x),\ x \in \Sigma _A$\ are positive, we prove the following fundamental theorems.

\medskip

\begin{thm}\label{thm1}
Suppose that $M$ is a H\"{o}lder continuous function on $\Sigma_A$ taking values in the set of positive $d\times d$ matrices. 
Then for any $q\in {\Bbb R}$, there is a unique $\sigma$-invariant, ergodic probability measure $\mu_q$ on $\Sigma_A$ of which one can find constants $C_1>0$, $C_2>0$ such that 
\begin{equation}\label{n1}
C_1\leq \frac{\mu_q([J])  }{\exp(-nP(q))\cdot\|M(x)M(\sigma x)\ldots M(\sigma^{n-1}x)\|^q}\leq C_2
\end{equation}
for any $n>0$, $J \in \Sigma_{A,n}$ and $x\in [J].$ 
\end{thm}
  
\medskip

The above measure $\mu_q$  is called the {\it Gibbs measure} associated with $M$ and $q$. We remark that the theorem generalizes the classical existence result of the Gibbs measure for a real-valued $M(x)$ (see [B, $\S$1.4]). The positivity of the matrices is used to  yield the follow simple estimate (Lemma 2.1)
\begin{eqnarray}
\|M(x)\cdots M(\sigma^{n+\ell-1}x)\|\approx
 \|M(x)\cdots M(\sigma^{n-1}x)\|\cdot\|M(\sigma^nx)\cdots
M(\sigma^{n+\ell-1}x)\|.
\end{eqnarray}
By using this we can apply a technique of Brown, Michon and J. Peyri\`ere [BMP] and Carleson [C] to construct a certain ergodic measure which is  the Gibbs measure $\mu_q$. The  $\mu_q$ has the following {\it quasi-Bernoulli property}  (Heurteaux [H]): there exists $C>0$ such that for any $ \ n, k \in {\Bbb N}$ with $ I\in \Sigma_{A,n} , \ J \in \Sigma_{A,k} $ and $ IJ \in \Sigma_{A, n+k}$
\begin{eqnarray}
C^{-1} \mu_q ([I])\mu_q([J]) \leq \mu_q([IJ]) \leq  C \mu_q ([I])\mu_q([J]) \qquad 
\end{eqnarray}
This together with a result of Heurteaux [H] imply

\medskip

\begin{thm}\label{thm2}  Under the condition of Theorem \ref{thm1}, $P(q)$ is differentiable for $q \not = 0$.
\end{thm}

\medskip

As an application, we let
$$
E(\alpha):= \big\{x\in \Sigma_A: \ \lim_{n\rightarrow \infty}
\frac{ \log \|M(x)M(\sigma x)\cdots M(\sigma^{n-1}(x)\|}{n}=\alpha \big\}.
$$ 
We prove the following dimension formula
\medskip

\begin{thm}\label{thma}
Under the same assumption of Theorem \ref{thm1}, we have for any $
\alpha=P'(q), \ q \not = 0$,
$$
\dim_HE(\alpha)=\frac{1}{\log m} (-\alpha q+P(q))
$$
where $\dim_H$ denotes the Hausdorff dimension.
\end{thm}

\medskip

The above theorems depend very much on the positivity of the matrix-valued $M(x)$. In order to extend them to nonnegative matrix-valued functions, we have to impose more conditions on $M(x)$:

(H1) \quad $M(x)=M_{i}$ \ if $x \in [i], \ i = 1, \cdots, m $;       

(H2) \quad $M$ is irreducible in the following sense: there exists $r > 0$ such that 

\quad \quad \quad for any  $i, j\in \{1,2,\ldots,m\}$,
\begin{equation}\label{n2}
\sum_{k=1}^r\sum_{K\in \Sigma_{A,k; i,j}}M_K > {\bf 0}
\end{equation}
\quad \quad \quad where $\Sigma_{A,k; i,j}$ denotes the set of all  $ K \in \Sigma_{A,k}$ such that $iKj\in \Sigma_{A,k+2}$.

\medskip

\noindent We see that under the assumption (H1), the pressure function in (1.1) can be re-written as 
\begin{equation}\label{n0}
P(q)=\lim_{n\rightarrow \infty} \frac{1}{n}\log 
\sum_{J\in \Sigma_{A,n}} \| M_J\|^q.
\end{equation}
where $M_J = M_{j_i} \cdots M_{j_n}$. If $\Sigma_A = \Sigma $ is the symbolic space with a full shift,  then (H2) is equivalent to $\sum_{k=1}^r (M_1 + \cdots+ M_m)^k > {\bf 0}$. 

\medskip

In this new setting, we use (H2) to adjust (1.3) and the required lemmas, the Gibbs measure $\mu_q$ is shown to  exist for $q >0$. This time $\mu_q$ only satisfies  $ \mu_q([IJ]) \leq C \mu_q([I])\mu_q([J])$ instead of (1.4); nevertheless we can still prove the differentiability of $P(q), \ q >0$ as in Theorem 1.2. Theorem 1.3 can be adjusted likewise (see Theorem 3.3 and Theorem 3.4). As an application, the first author proves the smoothness of the $L^q$-spectrum ($q>0$) and the multifractal formalism for a class of self-similar measures with overlaps (including the Bernoulli convolutions associated with Pisot numbers) in a forthcoming paper [Fe2].   

\bigskip
 For the organization of the paper, we prove the above results for the positive matrix-valued functions in Section 2. In Section 3, we modify the proofs for the non-negative matrix-valued functions with (H1) and (H2). In Section 4, we give an illustration of reducing an IFS with overlap to a vector-valued IFS with no overlap, and the pressure function in (1.6) arises. We also give some remark on the theorems and raise a few unsettled problems.

\bigskip

{\noindent}{\bf Acknowledgment}. The paper was originally written for the random products of $m$ matrices. The authors would like to thank the referee for the suggestion to modify it to the present form which can be appealed to more general situation. They also thank E. Olivier for introducing the multifractal results of \cite{BMP,H} and  \"{O}. Stenflo for reading the manuscript carefully and suggesting some improvements.

\section{\bf Positive Matrices}

In this section we assume that 
 $M$ is a H\"{o}lder continuous function on $\Sigma_A$ taking values in the set of all positive $d\times d$ matrices. 

For any two families of positive numbers  $\{a_i\}_{ i \in {\cal I}},\  \{b_i\}_{i \in {\cal I}}$, we write, for brevity,  $a_i \approx b_i$ to mean the existence of  a constant $C>0$ such that $C^{-1} a_i \leq  b_i  \leq C a_i$ for all $i \in {\cal I}$;
$a_i \preccurlyeq b_i$ to mean the existence of  a constant $C>0$ such that $ a_i \leq C b_i$ for all $i \in {\cal I}$ and  $a_i \succcurlyeq b_i$ means $b_i \preccurlyeq a_i$.

We start with a simple lemma:

\begin{lemma}\label{lem1}
For any $x\in \Sigma_A$, $n,\ell\in {\Bbb N}$,
\begin{eqnarray*}
 \|M(x)\cdots M(\sigma^{n+\ell-1}x)\|\approx
 \|M(x)\cdots M(\sigma^{n-1}x)\|\cdot\|M(\sigma^nx)\cdots
M(\sigma^{n+\ell-1}x)\|
\end{eqnarray*}
(the involved constant in $\approx$ is independent of $n, \ell$ and $x$).
\end{lemma}

\noindent \prf  \  It is clear that 
$$
\|M(x)\cdots M(\sigma^{n+\ell-1}x)\|
\leq \|M(x)\cdots M(\sigma^{n-1}x)\|\cdot\|M(\sigma^nx)\cdots
M(\sigma^{n+\ell-1}x)\|.
$$ 
To prove the reverse inequality, we observe that  $M$ is positive and continuous, there is a constant $C>0$ such that 
$$\frac{\min_{ i,j}M_{i,j}(x)}{\max_{ i,j}M_{i,j}(x)}\geq C \qquad \forall \ x\in \Sigma_A.
$$
This implies that $
M(x)\geq \frac{C}{d}EM(x)$ ($A \geq B$ means that $A_{i,j}\geq B_{i,j}$ for each index $(i,j)$) where $E=(E_{i,j})_{1\leq i,j\leq d}$ is the matrix 
 whose entries are all equal to $1$. Let ${\bf 1}$  be the $d$-dimensional column vector
 each coordinate of which is $1$. Then
\begin{eqnarray*}
\|M(x)\cdots M(\sigma^{n+\ell-1}x)\|
&\geq& \| M(x)\cdots M(\sigma^{n-1}x)\frac{C}{d}EM(\sigma^nx)\cdots
M(\sigma^{n+\ell-1}x)\|\\
&=& \frac{C}{d}\| M(x)\cdots M(\sigma^{n-1}x){\bf 1}^\tau{\bf 1}M(\sigma^nx)\cdots
M(\sigma^{n+\ell-1}x)\|\\
&=&\frac{C}{d}\| M(x)\cdots M(\sigma^{n-1}x)\|\cdot\|M(\sigma^nx)\cdots
M(\sigma^{n+\ell-1}x)\|.
\end{eqnarray*}
\hfill $\Box$

\bigskip

We define    
\begin {eqnarray*}
s_n(I,q)=\sup_{x\in [I]}\|M(x)M(\sigma x)\ldots M(\sigma^{n-1}x)\|^q\qquad \forall I\in \Sigma_{A,n} 
\end{eqnarray*}
and 
\begin{eqnarray}
\quad s_n(q)=\sum_{I\in \Sigma_{A,n}}s_n(I,q).
\end{eqnarray}
\begin{lemma}\label{lem2}  For a fixed  $q \in {\Bbb R}$, 
\begin{eqnarray*}
s_n(I, q) \approx  \|M(x)M(\sigma x)\ldots M(\sigma^{n-1}x)\|^q   \qquad \forall\ I\in \Sigma_{A,n}, \ x\in [I].
\end{eqnarray*}
\end{lemma}

\noindent \prf \ For any $n\in {\Bbb N}$, define 
\begin{equation}\label{en1}
\eta_n=\sup \left \{\frac{M_{i,j}(x)}{M_{i,j}(y)}:\ I\in \Sigma_{A,n},\  x,y\in [I], \ 1\leq i,j\leq d \ \right  \}. 
\end{equation}
Since each $M_{i,j}$ is  positive and H\"{o}lder continuous, we have $|\log \eta_n | \leq C m^{-\alpha n}$ for some $C>0$ and $0 < \alpha < 1$. It follows easily that 
$\eta : = \prod_{n=1}^\infty \eta_n < \infty $
and hence for $x\in [I]$,
$$
\|M(x)\ldots M(\sigma^{n-1}x)\|^q \leq s_n(I, q) \leq \eta^{|q|}\|M(x)\ldots M(\sigma^{n-1}x)\|^q.
$$
\hfill $\Box$

\bigskip

We have assumed that $A$ is primitive, there is an integer $p > 0$ such that $A^p > {\bf 0}$.  This implies that for any $I\in \Sigma_{A,n},J\in \Sigma_{A,\ell}$, there exists $K\in \Sigma_{A,p}$ such that $IKJ\in \Sigma_{A,n+\ell+p}$. 

\medskip

\begin{lemma}\label{lem3} Let $p$ be such that $A^p > {\bf 0}$. Then for a fixed  $q \in {\Bbb R}$, 

{} \ \ {\rm (i)} $s_\ell(q)\approx s_{\ell-p}(q)$ for all  $ \ell > p$;

{} \ {\rm (ii)} For  $ I\in\Sigma_{A,n}, \ \ell> p$,
$\sum_J s_{n+\ell}(IJ,q) \approx \sum_J s_{n+\ell}(JI,q) \approx s_n(I,q)s_\ell(q)$ \\
where the first (second) sum is taken over all $J\in \Sigma_{A,\ell}$ such that $IJ\in \Sigma_{A, n+\ell}$ \ ($JI\in \Sigma_{A, n+\ell}$ respectively ); 

{\rm (iii)}
$\sum_{K:\ IKJ\in \Sigma_{A,i}} s_i(IKJ,q) \approx s_n(I,q)s_\ell(J,q)s_{i-n-\ell}(q)$  
for all $I \in \Sigma_{A,n}, \ J\in\Sigma_{A,\ell},\ i>n+\ell+2p$.\\

\end{lemma}

\noindent \prf \ For any $I\in \Sigma_{A,\ell}$, write $I=KJ$ where $J\in \Sigma_{A,\ell-p}$.  By Lemmas \ref{lem1}, \ref{lem2},
we have (note that $p$ is fixed)
$$
s_{\ell}(I,q)\approx s_{\ell-p}(J,q).
$$
Since $A^p > {\bf 0}$, for $J\in \Sigma_{A,\ell-p}$, we can find $K \in \Sigma_{A, p}$ such that $I =KJ \in \Sigma_{A,\ell}$. Hence  when we take the sum  of $I \in\Sigma_{A,\ell}$ on the left side of the expression, it is \ $\approx$ \ to the right side  summing over all $J \in  \Sigma_{A,\ell-p}$. This implies (i).

To prove (ii), we fix $I\in \Sigma_{A,n}$ and  take $J\in \Sigma_{A,\ell}$ such that $IJ\in \Sigma_{A, n+\ell}$. By Lemmas \ref{lem2}, \ref{lem1}, we have 
$$
s_{n+\ell}(IJ,q)\approx s_n(I,q)s_\ell(J,q).
$$
Thus 
$$
\sum_J s_{n+\ell}(IJ,q)\preccurlyeq s_n(I,q)s_\ell(q).
$$
For the reverse inequality we note for any $J^\prime\in \Sigma_{A,\ell-p}$, there is $K\in \Sigma_{A,p}$ such that $IKJ^\prime\in \Sigma_{A,n+\ell}$ and 
$$
s_{n+\ell}(IKJ^\prime)\approx s_n(I,q)s_p(K,q)s_{\ell-p}(J^\prime,q)\approx s_n(I,q)s_{\ell-p}(J^\prime,q).
$$
Therefore summing over the above $J^\prime$, we have 
 $$
\sum_{IJ\in \Sigma_{A, n+\ell}}  s_{n+\ell}(IJ,q)\succcurlyeq \sum_{J^\prime} s_{n+\ell} (IKJ^\prime, q) \approx s_n(I,q)s_{\ell-p}(q)\approx s_n(I,q)s_{\ell}(q)
$$
(we make used of $A^p > {\bf 0}$ as in (i)).    This proves one of the \ $\approx $ \ in (ii). The remaining part follows from the same argument.

To prove (iii), we first observe that
\begin{eqnarray*}
\sum_{K:\ IKJ\in \Sigma_{A,i}} s_{i}(IKJ,q) 
& \approx & \sum_{K:\ IKJ\in \Sigma_{A,i}} s_{n}(I,q)s_{i-n-\ell}(K,q)s_{\ell}(J,q)\\
& \preccurlyeq & s_{n}(I,q)s_{\ell}(J,q)s_{i-n-\ell}(q)
\end{eqnarray*}
On the other hand, for any $K_1\in \Sigma_{A,i-n-\ell-2p}$, there exists $K_2,K_3\in \Sigma_{A,p}$ such that 
$IK_2K_1K_3J\in \Sigma_{A,i}$. Therefore  
\begin{eqnarray*}
\sum_{K:\ IKJ\in \Sigma_{A,i}} s_{i}(IKJ,q)&\succcurlyeq&
s_{n}(I,q)s_{\ell}(J,q)\sum_{K_1 \in \Sigma_{A,i-n-\ell-2p}}s_{i-n-\ell-2p}(K_1,q)\\
&\approx & s_{n}(I,q)s_{\ell}(J,q)s_{i-n-\ell-2p}(q)\\
&\approx & s_{n}(I,q)s_{\ell}(J,q)s_{i-n-\ell}(q).  
\end{eqnarray*}  
\hfill $\Box$

\medskip

\begin{lemma} \label{lem4} \ For a fixed $q \in {\Bbb R}$,   

{} \ {\rm (i)} $s_{\ell + n}(q) \approx s_\ell (q)s_n(q)$.

 {\rm (ii)} $s_n(q) \approx  \exp(nP(q))$ where $P(q)$ is the pressure function defined in (1.1).
\end{lemma}

\medskip

\noindent {\bf Proof.}\  From Lemma \ref{lem3} (ii),  there exist $C, C' >0$ such that  
$$
C' s_\ell (q)s_n (q) \leq s_{\ell +n}(q) \leq C s_\ell (q)s_n(q),
$$
which proves (i). To prove (ii),  we can write $C s_{\ell + n}(q) \leq  \left( C s_\ell (q)\right) \left(C s_n (q)\right)$. Hence the subadditivity property implies
$$
 P(q) =   \lim_{n\rightarrow \infty} \frac{ \log \left(C s_n (q)\right)}{n}= \inf_n \frac{ \log \left(Cs_n (q)\right)}{n},
$$
so that $ C^{-1}\exp(nP(q))\leq s_n(q).$
The reverse inequality follows from a similar argument.
\hfill   $\Box$

\bigskip

For each integer $n>0$, let ${\cal B}_n$ be the $\sigma$-algebra generated by the cylinders  $[I]$, $I\in \Sigma_{A,n}$. We define a sequence of probability measures $\{\nu_{n,q}\}$ on  ${\cal B}_n$ by 
\begin{equation}
\nu_{n,q}([I])=\frac{s_n(I,q)}{s_n(q)} \qquad \forall \ I\in \Sigma_{A,n}\ .
\end{equation}
 Then there is a subsequence $\{\nu_{n_k,q}\}_{k\geq 1}$ converges in the weak-star 
 topology to a probability measure $\nu_q$. The following assertion shows that $\nu_q$ has the ``Gibbs property''. 

\medskip
 
\begin{lemma}\label{lem5} For a fixed $q\in {\Bbb R}$,  $\nu_q([I]) \approx s_n(I,q)\exp(-nP(q))$ for all  $ n>0, \ I\in \Sigma_{A,n} .$
\end{lemma}

\noindent {\bf Proof.} \ Let $p$ be such that $A^p > {\bf 0}$. For any $I\in \Sigma_{A,n}$ and $\ell > n+p$, we have 
\begin{eqnarray*}
 \nu_{\ell,q}([I]) &=& \sum_{J:\ IJ\in \Sigma_{A,\ell}}\nu_{\ell,q}([IJ])
 =\sum_{J: \ IJ\in \Sigma_{A,\ell}} \frac{s_\ell(IJ,q)}{s_\ell (q)}\\
 &\approx & s_n(I,q)\frac{s_{\ell-n}(q)}{s_\ell (q)}  \qquad \qquad \qquad \qquad \hbox {(by  Lemma \ref{lem3} (ii))}\\
&\approx&  s_n(I,q)\exp(-nP(q)). \qquad   \qquad \quad \hbox {(by Lemma \ref{lem4})}.
 \end{eqnarray*}
Letting $\ell  = n_k \uparrow  \infty$, we obtain the desired result.
 \hfill   $\Box$
 
 \bigskip

\noindent {\bf Proof of Theorem \ref{thm1}.} \ 
Fix $q\in {\Bbb R}$.    Let  $\mu_q$ be  a limit point 
 of  the subsequence of $\left\{\frac{1}{n}(\nu_q+\nu_q\circ \sigma^{-1}
 +\ldots +\nu_q\circ \sigma^{-(n-1)})\right\}$ in the weak-star
 topology. Then $\mu_q$ is a $\sigma$-invariant measure on 
 $\Sigma_A$.  We have for each $I\in \Sigma_{A,n}$ and $\ell >p$, 
\begin{eqnarray}
  \nu_q\circ \sigma^{-\ell}([I]) & = & \sum _{J: \ JI\in \Sigma_{A,n+\ell}}\nu_q([JI]) \nonumber\\
 &\approx & \sum _{J:\ JI\in \Sigma_{A,n+\ell}}s_{n+\ell}(JI,q)\exp(-(n+\ell)P(q))\quad \hbox{( by Lemma \ref{lem5})} \nonumber\\
& \approx & s_\ell(q)s_n(I,q)\exp(-(\ell + n)P(q)) \qquad \qquad \hbox {(by  Lemma \ref{lem3} (ii))} \nonumber \\
 &\approx &  s_n(I,q)\exp(-nP(q)) \qquad \qquad  \qquad \qquad \hbox {(by Lemma \ref{lem4})}. 
\label{ee3}
\end{eqnarray}
This  proves that  $\mu_q$ is a Gibbs measure.   In what follows we prove that $\mu_q$ is ergodic. First we show that  there is a constant $C>0$ such that for each  $I\in \Sigma_n$, $J\in \Sigma_\ell$, 
\begin{equation}
\label{uuu}
\lim_{k\to \infty}\frac{1}{k}\sum_{i=0}^{k-1}\mu_q\Big([I]\cap\sigma^{-i}([J])\Big)\geq C\mu_q([I])\mu_q([J]). 
\end{equation}
Since $\mu_q$ is supported on $\Sigma_A$, it suffices to prove (\ref{uuu}) for $I\in \Sigma_{A,n}$ and $J\in \Sigma_{A,\ell}$. Note that when $i>n+2p$,
\begin{eqnarray*}
&&\mu_q\Big([I]\cap\sigma^{-i}([J])\Big)\\
&=&\sum_{K: \ IKJ \in \Sigma_{A,i+\ell}}  \mu_q([IKJ]) \\
&\succcurlyeq&
\sum_{K: \ IKJ \in \Sigma_{A,i+\ell}}s_{i+\ell}(IKJ,q)\exp(-(i+\ell)P(q))  \\
&\approx& s_n(I,q)s_{\ell}(J,q)s_{i-n}(q)\exp(-(i+\ell)P(q)) \qquad \quad \hbox { (by Lemma \ref{lem3} (iii))}\\
&\approx& s_n(I,q)s_{\ell}(J,q)\exp(-(n+\ell)P(q))\\
&\approx & \mu_q([I])\mu_q([J]) 
\end{eqnarray*}
from which (\ref{uuu}) follows.   Since the collection 
$\{[I]:I\in \Sigma_n,\ n\in {\Bbb N}\}$ is a semi-algebra that generates the Borel $\sigma$-algebra on $\Sigma$, a standard argument (e.g., see the proof of  [W, Theorem 1.17])  shows that 
for any Borel sets $A, B\subset \Sigma$, 
$$
\lim_{k \to \infty}\frac{1}{k}\sum_{i=0}^{k-1}\mu_q\Big(A\cap\sigma^{-i}(B)\Big)\geq C\mu_q(A)\mu_q(B). 
$$ 
This implies that for any Borel sets  $A, B\subset \Sigma$ with $\mu_q(A)>0$, $\mu_q(B)>0$, there exists $n>0$ with 
$\mu_q(A\cap \sigma^{-n}(B))>0$. By [W, Theorem 1.5], $\mu_q$ is ergodic.

For the uniqueness we recall that any two distinct ergodic measures must be singular to each other; but the Gibbs property (1.2) implies that any two $\mu_q$  must be absolutely continuous to each other. Hence $\mu_q$ must be unique.
\hfill $\Box$

\bigskip

\begin{cor}  Let $\mu_q$ be the Gibbs measure in Theorem 1.1. There exists $C>0$ such that for any $  I \in \Sigma_{A, n}, J \in \Sigma_{A, \ell}$ with $IJ \in \Sigma_{A, n+\ell}$, 
$$
C^{-1} \mu_q ([IJ]) \leq \mu_q([I]) \mu_q([J]) \leq C \mu_q ([I]) \mu_q([J]) .
$$
\end{cor}

\medskip

\noindent \prf  \   We have seen from the proof of Lemma 2.3 that for the above $I, J$, $s_{n+\ell}(IJ, q) \approx s_n(I,q) s_\ell (J, q)$ and from Lemma 2.4,
$s_{\ell +n}(q) \approx s_n(q) s_\ell(q)$. By the definition of $\nu_{n,q}$,  we have 
$$
\nu_{n,q} ([IJ]) \approx \nu_{n,q} ([I]) \nu_{n,q}([J])
$$
which implies that the Gibbs measure $\mu_q$ has the same property. 
\hfill $\Box$
\bigskip

The above property is called {\it quasi-Bernoulli property} by Heurteaux [H] (we remark that Heurteaux only introduced and studied it for measures in the full shift space $\Sigma$). To prove Theorem \ref{thm2}, we need a result in [H]. Let $\eta$ be a   probability measure on $\Sigma$.
For $q \in {\Bbb R}$, let $\tau_\eta(q)$ be the $L^q$-spectrum of $\eta$, i.e.,
$$
\tau_\eta(q)=\liminf_{n\rightarrow\infty} 
\frac{\log \sum_{I}\eta([I])^q }{\log m^{-n}},
$$
where the summation is taken over all $I\in \Sigma_n$ with $\eta([I])\neq 0$.

\medskip

 
\begin{pro} \label{H}{\rm (\cite[Theorem 2.1]{H})}
Let $\eta$ be a probability measure on $\Sigma$. Assume that there exists a 
constant $C>0$ such that 
\begin{equation}\label{22}
\eta([IJ])\leq C\eta([I])\eta([J])\qquad \forall \ I\in \Sigma_n,\ J\in \Sigma_\ell \ . 
\end{equation}
Then $\tau^\prime_\eta(1)$ exists if  $\eta$ is a Young measure 
(i.e., $\displaystyle \lim_{n\rightarrow \infty}\frac{\log \eta\big( I_n(x)\big)}{\log m^{-n}}={\rm constant}$ for $\eta$ almost all $x = (j_i) \in \Sigma$, here $I_n(x) =[j_1 \ldots j_n] )$. 
\end{pro}

\bigskip
\noindent {\bf Proof of Theorem \ref{thm2}.} \ 
For each $q$, let $\mu_q$ be the corresponding Gibbs measure in Theorem \ref{thm1}. We can view $\mu_q$ as a measure on $\Sigma$. For $t\in {\Bbb R}$, let  $\tau_{\mu_q}(t)$ be the $L^t$-spectrum of ${\mu_q}$. Since $\mu_q$ has the Gibbs property, it is easy to show by the definition of  $L^t$-spectrum that
\begin{equation}\label{ett}
\tau_{\mu_q}(t)=\frac{tP(q)-P(qt)}{\log m}.
\end{equation}
Note that $\mu_q$ satisfies the condition (\ref{22}).
Since $\mu_q$ is ergodic on $\Sigma$, it is a Young measure by the Shannon-McMillan-Brieman theorem (i.e., $\displaystyle \lim_{n\rightarrow \infty}\frac{- \log \mu_q \big( I_n(x)\big)}{n}$ equals the  entropy of $\mu_q$ (with respect to $\sigma$) for $\mu_q$-almost all $x = (j_i) \in \Sigma$ and $I_n(x) =[j_1 \ldots j_n] )$. Hence by Proposition 2.7, $\tau_{\mu_q}(t)$ is differentiable at $t=1$. This implies that $P(q)$ is differential at any fixed $q\neq 0$, and 
$$
P^\prime(q)=\frac{P(q)-\log m\cdot\tau_{\mu_q}^\prime(1)}{q}.
$$
\hfill $\Box$   


\bigskip

\noindent {\bf Proof of Theorem \ref{thma}.} \ Let $\alpha = \tau^\prime (q)$ with $q \not =0$. Let 
$\mu_q$ be the corresponding Gibbs measures in Theorem \ref{thm1}, then (2.7) implies that  
$$
\tau_{\mu_1}(q)=\frac{qP(1) -P(q)}{\log m} \qquad  
$$
and
$$
E(\alpha)=
\left \{x\in \Sigma:\ \lim_{n\rightarrow \infty} \frac{\log \mu_1([x_1\cdots x_n])}{
\log m^{-n}} =\frac{P(1)-\alpha}{\log m}\right \} 
$$
By [BMP, Theorem 1] or [LN2, Theorem 4.1], we have 
\begin{equation}\label{en}
\dim_HE(\alpha)\leq \big (\frac{ P(1)-\alpha}{\log m}\big )q - \tau_{\mu_1}(q)=\frac{1}{\log m}(-\alpha q+P(q))   \qquad    \forall \  q \in {\Bbb R}.
\end{equation}

 For the reverse inequality, we see from the proof of Theorem \ref{thm2} that $\tau_{\mu_q}'(1)$ exists   and 
$$
\tau_{\mu_q}^\prime(1)= \frac{-qP^\prime(q)+P(q)}{\log m}.
$$
By [N], we have for $ \mu_q $ almost all $I=(i_1\ldots i_n\ldots) \in \Sigma $,
$$
\lim_{n\to \infty}\frac{\log \mu_q([ i_1, \ldots  , i_n])}{\log m^{-n}}
= \tau^\prime_{\mu_q}(1) = \frac{-qP^\prime(q)+P(q)}{\log m}.
$$
This implies that 
$$
\lim_{n\to \infty}\frac{\log \|M(x)M(\sigma x)\cdots M(\sigma^{n-1}x)\|}{n}
=P^\prime(q) = \alpha \qquad \mu_q-a.a.\ x \in \Sigma.
$$
Therefore we have 
$$
\dim_HE(\alpha)\geq \dim_H \mu_q =
\frac{-q\alpha + P(q)}{\log m}.
$$
\hfill  $\Box$

\bigskip

\section{\bf Nonnegative matrices}

In this section, we always assume that  $M$ is a function on $\Sigma_A$  taking  values in the set of all $d\times d$ non-negative matrices and satisfies (H1) and (H2) defined in Section 1. Let $q>0$ be fixed. Then $s_n(I, q)$ and $s_n(q)$ in (2.1) are reduced to 
$$
s_n(I,q)=\|M_I\|^q \quad \forall \ I \in \Sigma_{A,n} \qquad \hbox {and} \qquad 
s_n(q) = \sum_{I\in \Sigma_{A,n}} s_n(I,q). 
$$
For convenience, we let 
\begin{equation}\label{f3}
b=\min_{s,t\in \{1,2,\ldots m\}}\min_{1\leq i, j\leq d}\left (\sum_{k=1}^r\sum_{K\in \Sigma_{A,k,s,t}}M_K \right )_{i,j}  \ .
\end{equation}
Then  $b>0$ by (H2).

We will reformulate the three theorems in the previous section. The proofs are almost the same and for simplicity, we only point out the differences.  Here Lemmas \ref{lem1}, \ref{lem2} do not hold anymore; on the other hand we can use (H2) to replace these lemmas to obtain an analog of Lemma \ref {lem3}:  

\medskip

\begin{lemma}\  For a fixed $q>0$, 

{} \ \ {\rm (i)} $s_{\ell+1}(q)\approx s_\ell(q)$.

{} \ {\rm (ii)} For  $I\in\Sigma_{A,n}$,
$\sum_J s_{n+\ell}(IJ,q) \approx \sum_J s_{n+\ell}(JI,q) \approx  s_n(I,q)s_\ell(q)$    where the first (second) sum is taken over all $J\in \Sigma_{A,\ell}$ such that $IJ\in \Sigma_{A, n+\ell}$ \ ( $JI \in \Sigma_{A, n+\ell}$ respectively).

{\rm (iii)}
$\sum_{k=1}^{2r}\sum_{K:\ IKJ\in \Sigma_{A,i+k}} s_{i+k}(IKJ,q)\approx s_n(I,q)s_\ell(J,q)s_{i-n-\ell}(q)$ 
 \ for all \ $I\in\Sigma_{A,n},\  J\in\Sigma_{A,\ell},\ i>n+\ell.$
 
\end{lemma}

\noindent {\bf Proof.}\ \ For any $I\in \Sigma_{A,\ell+1}$, write $I=iJ$ with $J\in \Sigma_{A,\ell}$.
Using $\|M_I\|\leq \|M_i\|\|M_J\|$, we have for $q>0$,
$$
s_{\ell+1}(q)\leq m \big (\sup_{i\in \{1,2,\ldots,m\}}\|M_i\|^q \big ) s_\ell(q).
$$
That is,
$s_{\ell+1}(q)\preccurlyeq s_\ell(q)$.
For the reverse inequality,  since  for any $J\in \Sigma_{A,\ell}$, 
$$
\sum_{k=1}^r\sum_{K\in \Sigma_{A,k}:\ KJ\in \Sigma_{A,\ell+k}}\|M_{KJ}\|  
=\big \|\big (\sum_{k=1}^r\sum_{K\in \Sigma_{A,k}:\ KJ\in \Sigma_{A,\ell+k}}M_K\big )M_J\big \|\geq b\|M_J\|,
$$
it follows that  
$$ 
\sum_{k=1}^r\sum_{K\in \Sigma_{A,k}:\ KJ\in \Sigma_{A,\ell+k}}\|M_{KJ}\|^q
\geq 
\big(\frac{b}{\sum_{k=1}^rm^k}\big)^q\|M_J\|^q.
$$
This combines with  $s_{\ell+1}(q)\preccurlyeq s_\ell(q)$ imply that  
$$
s_{\ell+1}(q) \succcurlyeq \sum_{k=1}^r s_{\ell+k}(q)\geq s_\ell(q)
$$
and completes the proof of (i).

To prove (ii), it follows from $\|M_{IJ}\|\leq \|M_I\|\|M_J\|$ that 
$\sum_{J\in \Sigma_{A,\ell}: IJ\in \Sigma_{A,n+\ell}}s_{n+\ell}(IJ,q)\\
\preccurlyeq s_n(I,q)s_\ell(q).$ For the reverse inequality, we  use (H2) as above to conclude that for any $J\in \Sigma_{A,\ell}$ ,  
$$
\sum_{k=1}^r\sum_{K\in \Sigma_{A,k}: IKJ \in\Sigma_{A,n+\ell+k}}\|M_{IKJ}\|\succcurlyeq \|M_I\|\|M_J\|.
$$
Hence
$ 
\sum_{k=1}^r\sum_{K\in \Sigma_{A,k}: IKJ \in\Sigma_{A,n+\ell+k}}\|M_{IKJ}\|^q\succcurlyeq \|M_I\|^q\|M_J\|^q
$
and therefore, summing over the $J\in \Sigma_{A,\ell}$,
$$
\sum_{k=1}^r\sum_{L\in\sum_{A,\ell+k}:IL\in \sum_{A,n+\ell+k}}s_{n+\ell+k}(IL,q)\succcurlyeq  s_n(I,q)s_\ell(q) .
$$
 Since 
$$
\sum_{J\in\sum_{A,\ell}:IJ\in \sum_{A,n+\ell}} s_{n+\ell}(IJ,q)\succcurlyeq  
\sum_{J\in\sum_{A,\ell+1}:IJ\in \sum_{A,n+\ell+1}} s_{n+\ell+1}(IJ,q) ,
$$
we have 
$$
\sum_{J\in\sum_{A,\ell}:IJ\in \sum_{A,n+\ell}}s_{n+\ell}(IJ,q)\succcurlyeq 
\sum_{k=1}^r\sum_{J\in\sum_{A,\ell+k}:IJ\in \sum_{A,n+\ell+k}}s_{n+\ell+k}(IJ,q)
\succcurlyeq  s_n(I,q)s_\ell(q).
$$
This completes the proof of an $\approx$ in (ii); the other \ $\approx$ \  follows from an identical argument. 

To prove (iii), we have,  
\begin{eqnarray*}
\sum_{k=1}^{2r} \sum_{K:IKJ\in \Sigma_{A,i+k}} s_{i+k}(IKJ,q)
& \leq & s_n(I,q)s_\ell(J,q)\sum_{k=1}^{2r}s_{i+k-n-\ell}(q) \\
& \approx & s_n(q)s_\ell(J,q)s_{i-n-\ell}(q)   \qquad \qquad (\hbox {by}\ (i)).
\end{eqnarray*}
On the other hand, for any $W\in \Sigma_{A, i-n-\ell}$, by (H2), there exist $1\leq k_1\leq r$, $K_1\in \Sigma_{A,k_1}$ such that $IK_1W\in \Sigma_{A,i-\ell+k_1}$ and 
$$
\|IK_1W\|\geq \frac{b\|I\|\|W\|}{\sum_{k=1}^rm^k},
$$
 where $b$ is defined by (\ref{f3}). By using (H2) again,    
there exist $1\leq k_2\leq r$,  $K_2\in \Sigma_{A,k_2}$ such that $IK_1WK_2J\in \Sigma_{A,i+k_1+k_2}$ and 
$$
\|IK_1WK_2J\|\geq \frac{b\|IK_1W\|\|J\|}{\sum_{k=1}^rm^k}\geq
\frac{b^2\|I\|\|W\|\|J\|}{(\sum_{k=1}^rm^k)^2}.
$$
Therefore we have 
$$
\sum_{k=1}^{2r}
\sum_{K:IKJ\in \Sigma_{A,i+k}} s_{i+k}(IKJ,q)
\succcurlyeq s_n(I, q)s_\ell(J,q)s_{i-n-\ell}(q)
$$ 
and (iii) follows.
\hfill $\Box$

\bigskip 

We now state the corresponding theorems as in Section 1.

\medskip

\begin{thm}\label{thm3} Suppose $M$ is a function on $\Sigma_A$  taking  values in the set of all $d\times d$ non-negative matrices and satisfies (H1) and (H2). Then for any $q>0$, 
there is a unique Gibbs measure $\mu_q$ on $\Sigma_A$ as in Theorem \ref {thm1}.
\end{thm}

\medskip

The proof is almost identical with that of Theorem \ref {thm1} using Lemma 3.1. The only adjustment  is to replace 
$$
\mu_q\left([I]\cap \sigma^{-i}([J])\right)\approx \mu_q(I)\mu_q(J),  \qquad i>n+2p. 
$$
 by 
$$
\sum_{i=1}^{2r}\mu_q\left([I]\cap \sigma^{-i}([J])\right)\approx \mu_q(I)\mu_q(J)\qquad \forall \ i>n.
$$
We use the same proof as in Section 2 for the next two theorems.

\bigskip

\begin{thm}\label{thm4} Under the same conditions of Theorem \ref{thm3}, $P(q)$ is differentiable for any $q>0$
\end{thm}
 \medskip

\begin{thm}\label{thmb}
Under the same conditions  of Theorem \ref{thm3}, we have for any $
\alpha=P'(q)$, $q>0$,
$$
\dim_HE(\alpha)=\frac{1}{\log m} (-\alpha q+P(q))
$$
where $E(\alpha)  = \big \{ J=(j_i) \in \Sigma_A :  \ \lim_{n\to \infty}  {\log \|M_{j_1} \cdots M_{j_n} \|}/ {n} \ = \               \alpha \big \}$. 
\end{thm}

\medskip

To relate Theorem 3.4 to the classical random product of matrices,  we let $\{Y_n\}$ be the i.i.d. random variables that take values ${M_1, \ldots , M_m}$ , invertible matrices  and with uniform distribution, then $ \lim_{n\rightarrow \infty} \frac 1n \log \| Y_n \ldots Y_1\|= \lambda $ a.s. and $\lambda$ is called the upper Lyapunov exponent ([FK], [BL, Chapter 1]). In comparison with Theorem 3.4, we let $\Sigma_A = \Sigma$ be the space of full shift (i.e., all the entries of $A$ are $1$), then $P(0) = \log m$. The limit of the random variables corresponds to the case for $q=0$,  $\lambda = P'(0)$ and $\dim_H E(\lambda) = P(0)/\log m =1$ (the existence of the derivative follows from some additional assumptions on the $M_j$ ([BL, p.119]).

\medskip

We remark that if $\Sigma_A = \Sigma$, then condition (H2) is reduced to a more simple form: $\sum_{k=1}^r H^k >{\bf 0}$ where $ H = M_1 + \cdots +M_m $. The condition is essential for the theorems in Section 3. Indeed we have

\bigskip  

\noindent {\bf  Example 3.5.}  Let $
M_1=\left(\begin{array}{ll}
2 & 0 \\
0 & 1
\end{array}\right) $
and $ M_2=\left(\begin{array}{ll}
2 & 0 \\
0 & 3
\end{array}\right). $
 Then $H = M_1 + M_2$ is reducible. Since $M_J = 
\left (\begin{array}{ll}
2^n & 0 \\
0 & 3^k
\end{array} \right ) $
where $|J|=n $ and $k$ is the number of $2$'s appeared in  $J$.  $\sum_{|J|=n} \|M_J\|^q = \sum_{k=0}^n\binom{n}{k}(2^n+3^k)^q$. 
Note that 
$$
\sum_{k=0}^n\binom{n}{k}(2^n+3^k)^q \geq \max \Big \{ \sum_{k=0}^n\binom{n}{k}2^{nq},\ 
\sum_{k=0}^n\binom{n}{k}3^{kq}\Big \}=
\max\{ 2^{n(q+1) },\ 
(1+3^{q})^n\}
$$
and 
$$
\sum_{k=0}^n\binom{n}{k}(2^n+3^k)^q\leq  \sum_{k=0}^n\binom{n}{k}2^q(2^{nq}+3^{kq})=2^q \big(2^{n(q+1)}+(1+3^{q})^n\big).
$$
We have  $P(q)=\max\{(q+1)\log 2, \ \log (1+3^q)\}$, which is not differentiable at $q=1$.

\bigskip

We see that the Gibbs measure $\mu_q$ in Section 2 has the quasi-Bernoulli property. However for the case of  non-negative matrices, only  $\mu_q ([IJ]) \leq C \mu_q([I]) \mu_q([J])$, $ I \in \Sigma_{n}, \ J \in \Sigma_{\ell}$. The following example shows that the reverse inequality may not hold. 

\medskip

\noindent {\bf  Example 3.6.} Let $ M_1=\left(
\begin{array}{ll}
1 & 1 \\
0 & 1
\end{array}
\right)$ and $M_2$ an arbitrarily positive matrix, then \ $H = M_1 + M_2$ \ is an irreducible positive matrix.  Let $J = 1\ldots 1$ ($n$-times), then $\|M_J\| = n+2$ and hence 
$$ 
  \ \| M_J\| \ \|M_J\| \ \geq \  \frac n2 \|M_{JJ}\| \ .  
$$
Since $\mu_q (I) \approx \exp(-nP(q))\cdot \|M_I\|, \ I \in \Sigma_n$, we see that 
there does not exist $C'>0$ such that $ C' \mu_q ([J]) \mu ([J]) \leq \mu ([JJ])$.

\bigskip

\section{\bf Examples and remarks}\setcounter {equation}{0}

Consider the classical Bernoulli convolution $X = (1-\rho)^{-1}\sum_{n=1}^\infty \rho^n X_n$ where the $X_n$'s  are i.i.d. random variables which take values ${0, 1}$ and with probability $1/2$ on each value. Let $\mu_\rho$ be the distribution measure of $X$. It is well known that for  $0 < \rho < 1/2$, the measure is a Cantor type measures. It was proved recently that $\mu_\rho$ is absolutely continuous for almost all $1/2< \rho <1$ [S], however, it is still not clear which $\mu_\rho$ is absolutely continuous or singular. The question has been subjected to intensive investigation, the reader can refer to the survey articles [L], [PSS] and the references there.  The interest of the Bernoulli convolution in our setting is that the $\mu_\rho$ satisfies the {\it self-similar} identity
$$
\mu_\rho = \frac 12 \mu_\rho S_1^{-1} +  \frac 12 \mu_ \rho S_2^{-1}
$$
where $S_1 x = \rho x , \ S_2x = \rho x + (1-\rho) $; $\{S_1, S_2\}$  is the iterated function system (IFS). The support of $\mu_q$ is $[0, 1]$. For $0 < \rho < 1/2$, the $S_J(0,1)$'s are disjoint (as in the basic intervals of the Cantor set); for $1/2 < \rho < 1$, the $S_J(0,1)$ overlaps which is the source of difficulty. 

For $\rho = (\sqrt 5 -1)/2$, the reciprocal of the golden ratio, it was shown by Erd\"os that $\mu_\rho$ is singular. In order to consider the multifractal structure of $\mu_\rho$,  we can put the IFS $\{S_1, S_2\}$  to a new set of IFS  $\{R_i\}_{i=1}^3$with no overlap:
$$
R_1(x) = \rho^2 x, \quad  R_2(x) = \rho^3x + \rho^2, \quad \hbox {and} \quad 
R_3(x) = \rho^2 x + \rho.
$$
Then the measure $\mu_\rho$ satisfies 
$$
\mu_\rho ([i_1\ldots i_n]) \approx \frac 1 {4^n} \| M_{i_1}\ldots M_{i_n}\|
$$
 where $ [i_1\ldots i_n] = R_{i_1}\ldots R_{i_n}([0,1])$ and
\[
M_1=\left(
\begin{array}{ll}
1 & 1 \\
0 & 1
\end{array} \right),  \qquad 
M_2=\left(\begin{array}{ll}
\frac 12 & \frac 12 \\
\frac 12 & \frac 12
\end{array}
\right) \qquad \hbox{and}\qquad 
M_3=\left(\begin{array}{ll}
1 & 0 \\
1 & 1
\end{array} \right)\]
([LN1], [Fe1], [FeO]). The $\{M_i\}_{i=1}^3$ satisfies the conditions (H2). Hence by Theorem 3.2, 3.3, $P(q)$ is differentiable for $q >0$ and  the multifractal formalism holds. 

  Actually more can be said about the $L^q$-spectrum $\tau (q)$ of $\mu_q$: an explicit formula was given in [LN1] for $q>0$ and was extended to $q<0$ in [Fe1]. By using the formula it was proved that $\tau (q)$ is differentiable (actually real analytic) on $\Bbb R$ except one point in ${\Bbb R}^-$ in [Fe1]. 

\bigskip
 
The above example of Bernoulli convolution  gives rise to another interesting question.  Note that the above example is a special case of the overlapping IFS that can be reduced to new sets of IFS with no overlap and the calculation of the $\tau(q)$ can be converted into the product of matrices as in (1.4). Such IFS forms an important subclass of those that satisfy the {\it weak separation condition} ([LN2], [LNR])(it will be interesting to classify this subclass of IFS). Under the weak separation  condition  it was proved that the multifractal formalism is valid provided that $\tau (q)$ is differentiable [LN2].  However we do not know its differentiability in the general case. In a forthcoming paper [Fe2], the first author proves the differentiability of $\tau(q)$ for $q>0$ in the case that the IFS is equicontractive and satisfies the {\it finite type condition} (see [NW]).

\bigskip

The behavior for $q<0$ is also important for the multifractal analysis. There is no problem when $M$ is a {\it positive}
matrix-valued as we consider in Section 2. For the non-negative matrix-valued $M$, $M_J$ can be ${\bf 0}$, we have to modify the pressure function $P(q)$ in (1.6) slightly: 
\begin{eqnarray}
P(q) = 
\lim_{n \rightarrow  \infty} \frac 1n {\log \sum_{J \in {\mathcal N}_n} \| M_J \|^q}  
\end{eqnarray}
where ${\mathcal N}_n$ consists of all the $J \in \Sigma_{A,n}$ such that $M_J \not = {\bf 0}$. It is clear that if $M_J \not = {\bf 0}$ for all $J \in \Sigma_{A,n}$, then the super-additivity of the sum in (4.1)  implies that the limit exists. We include a simple proposition with $\Sigma_A = \Sigma$ to set up the consideration: 

\medskip

\noindent {\bf Proposition 4.4} \  Suppose $M_1, \cdots, M_m$ are non-negative matrices and $H = \sum_{i=1}^mM_i$ is irreducible, then 
the limit in (4.1) exists for each $q<0$ .

\medskip

\noindent {\bf Proof.} By the irreducibility, there exists integer $r$ with $\sum_{k=1}^rH^k>0$. Hence there is a constant $C > 0$ such that for any two finite indices $I,J $, there exists $K_0 \in \bigcup_{k=1}^r\Sigma_k$  satisfying 
\begin{equation}\label{equ1}
0 < \|M_{IK_0J}\|\leq C \|M_I\|\|M_J\|.
\end{equation}
Denote by $s_n=\sum_{J\in {\mathcal N}_n}\|M_J\|^q$.
Then (4.2) implies $s_ns_\ell\leq C^{-q}\sum_{k=1}^r
s_{n+\ell+k}$. From (4.2) we also deduce that  for any finite index $I$, there exists $i\in \Sigma_1$ such that $M_{Ii}\neq 0$; Since  $\|M_{Ii}\|\leq C_1\|M_I\|$ for some constant $C_1>0$,  we have $s_n\leq C_1^{-q}s_{n+1}$ for any integer $n,\ell$.
It follows that  $s_ns_\ell\leq C^\prime s_{n+\ell+r}$ for some constant $C^\prime 
> 0$ (depending on $q$), which implies that $a_n=\frac{1}{C^\prime} s_{n-r}$ is super-multiplicative. This yields the existence of the limit.
\qquad
\hfill $\Box$.

\medskip
The differentiability of such $P(q)$ for $q<0$ is still unknown. We know that in the above Bernoulli convolution of the golden ratio,  it is possible for the $P(q)$ to be  non-differentiable at a point of $q < 0$ [Fe1]. On the other hand, it is known that by imposing some stronger conditions on the matrices, the pressure function $P(q)$ is analytic near $q =0$  (see e.g., [BL, Theorem 4.3]).

\bigskip

Finally we remark that we do not know whether the theorems can be extended  matrices with entries in ${\Bbb R}$. An important theorem concerning this is in [BL, Theorem 4.3] for the analyticity of $P(q)$ near zero. Much closer to our development is the scaling functions: $f(x) = \sum _{i=0}^m c_i f(2x -i)$. It is known that such function can be put into matrix form as in the previous example [DL1].  Daubechies and Lagarias studied the  multifractal formalism of the well known scaling function $D_4$ [DL2]. They showed the differentiability of the corresponding $\tau (q)$, but the consideration depends on the two $2\times 2$  matrices involved to have a common eigenvector. There are some extensions in [LM].

\bigskip

 \vskip 0.2cm 
{\footnotesize
\noindent De-Jun FENG: Department of Mathematical Sciences,
Tsinghua University, Beijing, 100084, P.R. China.  \\
 and \\
Institute of Mathematical Sciences, The Chinese University of
Hong Kong,  Hong Kong\\
 {\it E-mail}: \  dfeng@@math.tsinghua.edu.cn
\vskip 12pt

\noindent Ka-Sing Lau: Department of Mathematics, the Chinese
University of Hong Kong,  Hong Kong\\
 {\it E-mail}: kslau@@math.cuhk.edu.hk \vskip 12pt }

\end{document}